\documentclass[pdflatex,sn-mathphys-num]{sn-jnl}

\usepackage{graphicx}%
\usepackage{amsmath,amssymb,amsfonts}%
\usepackage{amsthm}%
\usepackage{booktabs}%

\theoremstyle{thmstyleone}%
\newtheorem{theorem}{Theorem}[section]%
\newtheorem{lemma}[theorem]{Lemma}%
\newtheorem{proposition}[theorem]{Proposition}%

\theoremstyle{thmstylethree}%
\newtheorem{remark}[theorem]{Remark}%

\begin{document}

\title[A Tetration Distribution]{A Tetration Distribution}

\author*[1]{\fnm{Micah A.} \sur{Thornton}}\email{mthornton11@twu.edu}

\affil*[1]{\orgdiv{Department of Mathematics}, \orgname{Texas Woman's University}, \orgaddress{\city{Denton}, \state{Texas}, \country{USA}}}

\abstract{We introduce a one-parameter family of continuous distributions on the unit interval whose cumulative distribution functions are given by the odd-height power towers \(x\uparrow\uparrow n\) (\(n=1,3,5,\dots\)).
The associated densities admit a simple recursive formula obtained by logarithmic differentiation.
We prove that every such density is strictly positive on \((0,1)\) and integrates to one, thereby establishing that the construction yields a genuine probability distribution.
Numerical investigation reveals a characteristic peak--dip--rise shape: a sharp mode lying just below \(e^{-e}\) for all examined heights, followed by a trough and a gradual recovery toward density value~1 at the right endpoint.
We then identify the limiting distribution as the height tends to infinity through odd integers.
Its cumulative distribution function is \(W(t)/t\), with \(t=\ln(1/x)\) and \(W\) the Lambert \(W\) function, on \([e^{-e},1]\), and is given by the lower branch of the attracting two-cycle of \(y\mapsto x^{y}\) below \(e^{-e}\).
The constant \(e^{-e}\) is exactly Euler's convergence boundary for infinite power towers, and the period-doubling bifurcation occurring there produces a one-sided inverse-square-root singularity in the limiting density; this explains the observed spike and shows that the modes converge to \(e^{-e}\) from the left.
At the critical point the limit law satisfies \(F_\infty(e^{-e})=W(e)/e=e^{-1}\) exactly, so the limiting distribution places mass precisely \(1/e\) below \(e^{-e}\).
The limiting mean and variance are obtained as explicit integrals with numerical values \(0.3686906\dots\) and \(0.1112922\dots\); we note that the mean is close to, but provably distinct from, \(1/e\).
We then develop inference for the family: the height is a discrete parameter, so maximum likelihood is an exact profile search, and in a two-parameter exponentiated extension the second parameter is eliminated in closed form.
A Monte Carlo study confirms that the estimators are well behaved, but fitting two real unit-interval datasets shows the family to be decisively outperformed by the Beta and Kumaraswamy distributions.
We identify the structural reason: the density takes the value one at both endpoints for every height, so the family cannot adapt either tail, and we delimit the narrow regime in which it could be appropriate.
Qualitative morphological analogies with certain competence and hype-cycle curves are noted briefly, but formal statistical modelling of such phenomena is left for future work.}



\maketitle

\section{Introduction}
\label{sec:intro}

In considering the natural operators that we typically work with in mathematical
settings, we often overlook their iterative extensions---the so-called
hyper-operators. They extend the original operators in an iterative way.
Consider the following:
\begin{equation}
    S(x) = x + 1.
    \label{eq:successor}
\end{equation}
In Equation~\eqref{eq:successor} a few symbols are taken for granted, particularly
the binary function of addition represented by the `$+$' symbol, and the value of
the Arabic numeral `$1$'. The function $S(x)$ is known as the \textit{successor}
function, and it simply returns the argument of the function incremented by one.
\begin{equation}
    A(x,y) = x+y = \underbrace{S(S(\dots S(x)\dots))}_{y\text{ times}} = \underbrace{S(S(\dots S(y)\dots))}_{x\text{ times}}
    \label{eq:addition}
\end{equation}
In Equation~\eqref{eq:addition} the `$+$' symbol is defined, and we can see how
addition is simply the iterative application of succession.
Multiplication also extends naturally if we follow the sequence:
\begin{equation}
    M(x,y) = x \cdot y = \underbrace{A(x,A(x, \dots A(x,x) \dots ))}_{y\text{ times}} = \underbrace{A(y,A(y, \dots A(y,y) \dots ))}_{x\text{ times}}.
    \label{eq:multiplication}
\end{equation}
Next we come to exponentiation, which is not commutative, and
therefore the second equality from Equations~\eqref{eq:addition} and~\eqref{eq:multiplication}
does not hold. Exponentiation is \textit{right-associative}:
\begin{equation}
    E(x, y)= x^y = \underbrace{M(x,M(x, \dots M(x,x) \dots ))}_{y\text{ times}} \ne \underbrace{M(y,M(y, \dots M(y,y) \dots ))}_{x\text{ times}}.
    \label{eq:exponentiation}
\end{equation}
Finally we come to the operation of tetration, upon which the
distribution introduced in this paper is based:
\begin{equation}
    T(x,y) = {}^yx = x \uparrow \uparrow y = \underbrace{E(x,E(x, \dots E(x,x) \dots ))}_{y\text{ times}} \ne \underbrace{E(y,E(y, \dots E(y,y) \dots ))}_{x\text{ times}}.
    \label{eq:tetration}
\end{equation}
The name tetration comes from the root prefix `tetra-' and the suffix `-ation' and
refers to the nature of the operator as arising from iterating previous operators.
Tetration is the first of the hyper-operators, a hierarchy first introduced by
Goodstein~\cite{Goodstein1947}; the double up-arrow notation is due to
Knuth~\cite{Knuth1976}. The convergence theory of infinite power towers goes back
to Euler~\cite{Euler1783}, who showed that the iteration \(y\mapsto x^{y}\)
converges precisely for \(x\in[e^{-e},e^{1/e}]\); the constant \(e^{-e}\) will
play a central role in what follows. A detailed survey of iterated exponentials
and infinite power towers was given by Knoebel~\cite{Knoebel1981}. More recent
probabilistic investigations of random power towers appear in
Dalthorp~\cite{Dalthorp2024}. Ackermann's foundational 1928
work~\cite{Ackermann1928} provided an early functional form that captures
iterated operators. Independently, Dijoux~\cite{Dijoux2020} constructed a family
of lifetime distributions based on continuous iteration of the
exponential-minus-one function; the odd-height tetration CDF formulation of the
present paper is related to but distinct from that construction.

Here we introduce a distribution that arises when we consider the map
\begin{equation}
    x \mapsto x \uparrow \uparrow n, \qquad x \in [0,1], \qquad n = 2k+1, \;\; k \in \mathbb{Z}_{\ge 0}
    \label{eq:cdfmap}
\end{equation}
on the unit interval. For a random variable $X$, the above map defines a valid
cumulative distribution function, and its derivative therefore defines a valid
probability density function:
\begin{equation}
    P(X \le x \;;\; n) = x \uparrow\uparrow n, \qquad n = 2k+1, \;\; k \in \mathbb{Z}_{\ge 0}.
    \label{eq:cdfdef}
\end{equation}

This is the key insight of the paper. We derive a closed recursive form for the
derivative in Section~\ref{sec:cons} and prove that it is always positive on
$(0,1)$ and integrates to one. The resulting density exhibits a characteristic
shape---a steep initial rise to a peak near the origin, followed by a dip and a
gradual recovery on the unit interval---which is investigated in
Section~\ref{sec:dens}. Section~\ref{sec:limit} identifies the limiting
distribution as the height tends to infinity and shows that the location of the
spike is governed by Euler's constant \(e^{-e}\).
Section~\ref{sec:inference} develops estimation and simulation for the family,
reports a Monte Carlo study, fits two real datasets, and delimits the range of
problems for which the family is---and is not---appropriate.
Section~\ref{sec:apps} notes a few qualitative morphological analogies.
Section~\ref{sec:conc} summarizes the results and lists open questions.

\section{Construction and Definition}
\label{sec:cons}

\begin{figure}[h!]
    \centering
    \includegraphics[scale=0.55]{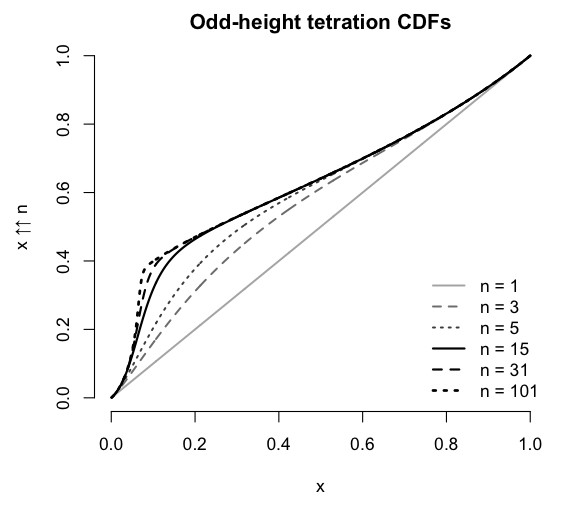}
    \caption{Cumulative distribution functions for selected odd-height tetration distributions.}
    \label{fig:cdf}
\end{figure}

Consider the form of the cumulative distribution function. Several members of
the family are plotted in Figure~\ref{fig:cdf}. These curves begin at zero and
increase monotonically to one over their domain, so they are candidates for
CDFs; their derivatives will be valid probability density functions provided
they remain positive on \((0,1)\). We establish this positivity below.

Throughout we write
\begin{equation}
    y_1(x) = x, \qquad y_m(x) = x^{\,y_{m-1}(x)} \quad (m \ge 2),
    \label{eq:towerdef}
\end{equation}
so that \(y_m(x) = x\uparrow\uparrow m\) for every positive integer \(m\),
\emph{whether or not \(m\) is odd}. The symbols \(F_n\) and \(f_n\) are reserved
exclusively for odd \(n\), where
\begin{equation}
    F_n(x) := y_n(x), \qquad f_n(x) := y_n'(x),
    \label{eq:Fdef}
\end{equation}
denote the CDF and density of the height-\(n\) tetration distribution. The
distinction matters: when \(m\) is even the function \(y_m\) is a perfectly
well-defined real-analytic function on \((0,1)\), but it is not a CDF, since it
is not non-decreasing on the whole of \([0,1]\) (indeed \(x^x \to 1\) as
\(x\to 0^{+}\)). We return to this point in Remark~\ref{rem:evenindex}. To prove
validity for odd heights it is enough to verify the Kolmogorov axioms on the
unit interval. We first derive the density recursion.

\begin{theorem}
\label{thm:recursion}
The derivative of the height-\(m\) power tower admits the recursive formula
\[
\frac{dy_m}{dx} = y_m\left( \frac{dy_{m-1}}{dx}\ln x + \frac{y_{m-1}}{x} \right)
\]
for \( m\ge 2 \), where \( y_1 = x \) and \( y_m = x^{y_{m-1}} \).
\end{theorem}

\begin{proof}
The relation \( y_m = x^{y_{m-1}} \) implies
\[
\ln y_m = y_{m-1}\ln x.
\]
Differentiating both sides with respect to \( x \) and multiplying through by \( y_m \) immediately yields the claimed recursion. The base case \( y_1 = x \) gives \( \frac{dy_1}{dx} = 1 \).
\end{proof}

Throughout the remainder of this section we fix \( x\in(0,1) \) and write
\( t = \ln(1/x) \in (0,\infty) \), so that the tower recursion becomes
\( y_1 = x \) and \( y_m = e^{-t\,y_{m-1}} \) for \( m\ge 2 \). The following
elementary bounds on the towers are the key to the positivity proof.

\begin{lemma}
\label{lem:evenbound}
For every \( x\in(0,1) \) and every integer \( k\ge 1 \),
\[
y_{2k}(x) > \frac{1}{e}
\qquad\text{and}\qquad
t\,y_{2k+1}(x) < 1 .
\]
\end{lemma}

\begin{proof}
We induct on \( k \). For the base case, \( y_2 = x^x \) attains its minimum on
\( (0,1) \) at \( x = e^{-1} \), where it equals \( e^{-1/e} > e^{-1} \); hence
\( y_2 > 1/e \).

Now suppose \( y_{2k} > 1/e \) for some \( k\ge 1 \). Then
\[
t\,y_{2k+1} = t\,e^{-t\,y_{2k}} < t\,e^{-t/e} \le \sup_{s>0}\, s\,e^{-s/e} = 1,
\]
where the first inequality is strict because \( y_{2k} > 1/e \), and the
supremum is attained at \( s=e \) with value \( e\cdot e^{-1} = 1 \).
Consequently
\[
y_{2k+2} = e^{-t\,y_{2k+1}} > e^{-1},
\]
which completes the induction.
\end{proof}

\begin{theorem}
\label{thm:positivity}
For every odd positive integer \( n \) and every \( x\in(0,1) \),
\[
\frac{d}{dx}\bigl(x\uparrow\uparrow n\bigr) > 0.
\]
\end{theorem}

\begin{proof}
Fix \( x\in(0,1) \) and define the scaled logarithmic derivatives
\[
g_m := \frac{x\,y_m'}{y_m}, \qquad m\ge 1 .
\]
Dividing the recursion of Theorem~\ref{thm:recursion} by \( y_m \) and
multiplying by \( x \) gives
\begin{equation}
g_m = y_{m-1}\bigl(1 + g_{m-1}\ln x\bigr) = y_{m-1}\bigl(1 - t\,g_{m-1}\bigr),
\qquad m\ge 2,
\label{eq:grec}
\end{equation}
with base case \( g_1 = 1 \). Since \( x>0 \) and \( y_m>0 \) on \( (0,1) \),
the derivative \( y_m' \) has the same sign as \( g_m \); it therefore suffices
to prove \( g_n>0 \) for every odd \( n \). We show by induction over the odd
integers that \( g_{n-2}>0 \) implies \( t\,g_{n-1}<1 \), which in turn gives
\( g_n = y_{n-1}(1-t\,g_{n-1}) > 0 \).

\medskip\noindent\emph{Base cases.} First, \( g_1 = 1 > 0 \). For \( n=3 \) we have
\( g_2 = y_1(1-t\,g_1) = x(1-t) \). If \( t\ge 1 \) then \( g_2\le 0 \), so
\( t\,g_2 \le 0 < 1 \); if \( 0<t<1 \) then \( t\,g_2 = t\,x(1-t) < t < 1 \).
In both cases \( t\,g_2<1 \), hence \( g_3 = y_2(1-t\,g_2) > 0 \).

\medskip\noindent\emph{Inductive step.} Let \( n\ge 5 \) be odd and assume \( g_{n-2}>0 \). We
distinguish two cases for the even-index quantity \( g_{n-1} \).
\begin{itemize}
\item[] \textbf{Case 1:} \( g_{n-2} \ge 1/t \). Then
\( g_{n-1} = y_{n-2}(1-t\,g_{n-2}) \le 0 \), and therefore
\( t\,g_{n-1} \le 0 < 1 \).
\item[] \textbf{Case 2:} \( 0 < g_{n-2} < 1/t \). Then \( 0 < 1-t\,g_{n-2} < 1 \),
so \( 0 < g_{n-1} < y_{n-2} \). Since \( n-2\ge 3 \) is odd,
Lemma~\ref{lem:evenbound} gives \( t\,y_{n-2} < 1 \), whence
\( t\,g_{n-1} < t\,y_{n-2} < 1 \).
\end{itemize}
In either case \( t\,g_{n-1}<1 \), and therefore
\( g_n = y_{n-1}(1-t\,g_{n-1}) > 0 \). By induction \( g_n>0 \), and hence
\( y_n' > 0 \), for every odd \( n \).
\end{proof}

\begin{remark}[Status of the even-index towers]
\label{rem:evenindex}
The recursion of Theorem~\ref{thm:recursion}, and equivalently the density
recursion~\eqref{eq:densityrec} below, evaluates the height-\(n\) density for odd
\(n\) in terms of the quantities \(y_{n-1}\) and \(y_{n-1}'\) of even index. It
is worth stating explicitly that this involves no circularity and no appeal to
any object of doubtful standing.

For even \(m\), the function \(y_m\) defined by~\eqref{eq:towerdef} is a
composition of finitely many real-analytic maps and is therefore itself
real-analytic, finite and strictly positive on \((0,1)\); likewise \(y_m'\) is
real-analytic and finite there. What fails for even \(m\) is \emph{only} the
monotonicity required of a CDF, and consequently \(y_m'\) is not a probability
density: by~\eqref{eq:grec}, \(g_2 = x(1-t)\) is negative exactly when \(t>1\),
that is on \(x\in(0,e^{-1})\), so \(y_2'<0\) there. The even-index towers are
used here purely as algebraic intermediates in a recursion among analytic
functions, and carry no probabilistic interpretation whatsoever. This is why the
notation of~\eqref{eq:Fdef} reserves \(F\) and \(f\) for odd indices.

Crucially, the proof of Theorem~\ref{thm:positivity} is arranged so that
\emph{no} positivity is ever required of an even-index quantity. Case~1 of the
inductive step explicitly admits \(g_{n-1}\le 0\), and the conclusion
\(g_n = y_{n-1}(1-t\,g_{n-1})>0\) is reached in that case precisely
\emph{because} \(g_{n-1}\) is non-positive, since then \(1-t\,g_{n-1}\ge 1>0\)
and \(y_{n-1}>0\). The sign changes of the even-index derivatives are thus
harmless: the induction passes through them rather than around them, and the
positivity of the odd-index densities is established independently of the sign
of any intermediate term. The same observation applies verbatim to the numerical
evaluation of~\eqref{eq:densityrec}, which involves only finite analytic
quantities at every step.
\end{remark}

\begin{theorem}
\label{thm:integrates}
For every odd positive integer \( n \),
\[
\int_0^1 \frac{d}{dx}\bigl(x\uparrow\uparrow n\bigr)\,dx = 1.
\]
\end{theorem}

\begin{proof}
Let \( F(x) = x\uparrow\uparrow n \).
By Theorem~\ref{thm:positivity}, \( F \) is continuously differentiable and
strictly increasing on \( (0,1) \). For odd \( n \) the tower extends
continuously to the closed interval with the boundary values
\[
F(0) = 0 \qquad\text{and}\qquad F(1) = 1.
\]
The fundamental theorem of calculus therefore yields
\[
\int_0^1 F'(x)\,dx = F(1) - F(0) = 1 - 0 = 1.
\]
\end{proof}

\section{Properties of the Density}
\label{sec:dens}

\begin{figure}[h!]
    \centering
    \includegraphics[width=\textwidth]{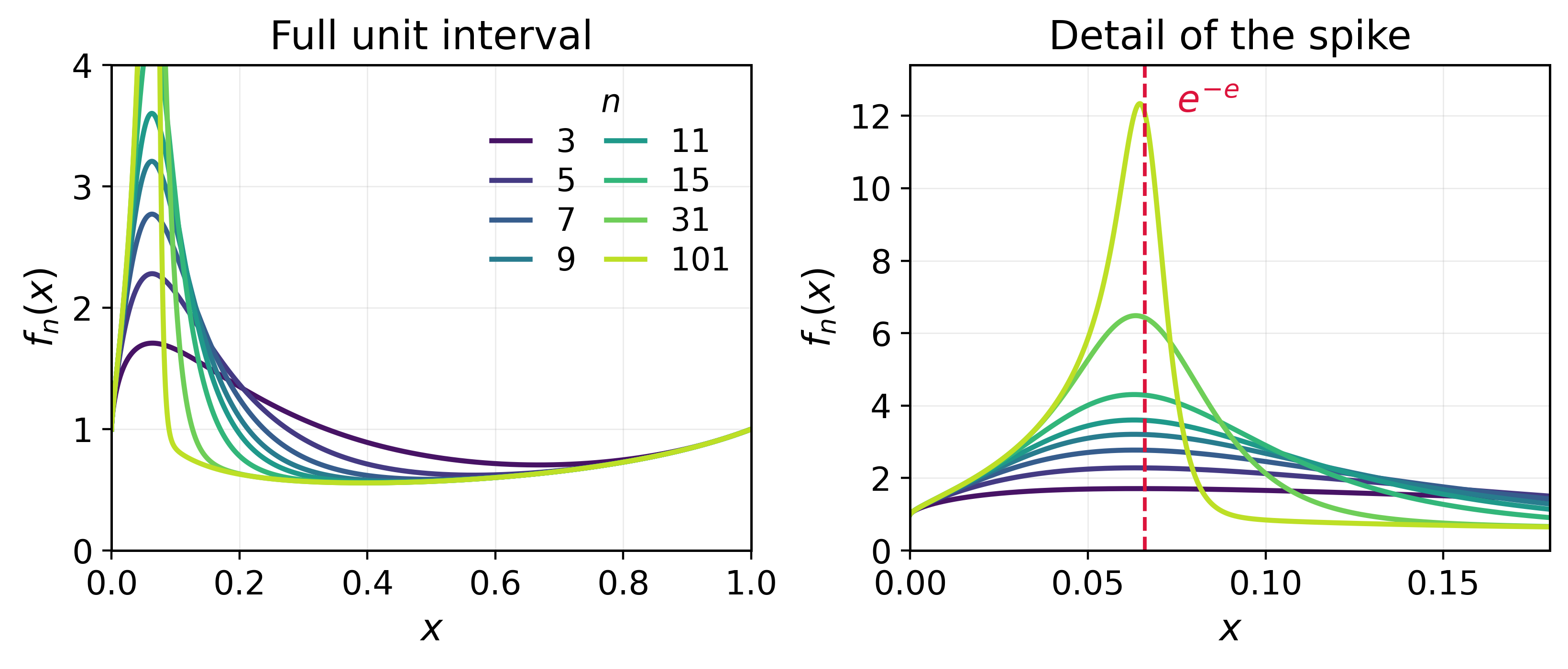}
    \caption{Probability density functions for selected odd-height tetration
    distributions. Left: the full unit interval. Right: detail of the spike,
    with the vertical dashed line marking Euler's constant
    \(e^{-e}\approx 0.065988\). As the height increases the peak grows taller and
    narrower and its location migrates toward \(e^{-e}\) from the left; see
    Section~\ref{sec:limit}.}
    \label{fig:pdf}
\end{figure}

Probability density functions for several of these distributions are shown in Figure~\ref{fig:pdf}.
The density \( f_n \) is obtained from the recursive formula derived in Section~\ref{sec:cons}:
\begin{equation}
    f_{n}(x) = F_{n}(x)\left( y_{n-1}'(x)\ln x + \frac{y_{n-1}(x)}{x}\right),
    \qquad f_1(x) = 1,
    \label{eq:densityrec}
\end{equation}
where \( F_n(x) = x\uparrow\uparrow n \) denotes the CDF and \(y_{n-1}\) is the
even-index auxiliary tower of~\eqref{eq:towerdef}; see
Remark~\ref{rem:evenindex}.
From the graphs one observes a characteristic shape: an initial spike near zero,
followed by a dip below height one, and then a gradual rise.
The height of the power tower primarily affects the sharpness of the initial spike
and how quickly the density falls below one.

\subsection{Mode}
\label{sec:mode}
Numerical maximization of the density yields the approximate modes reported in
Table~\ref{tab:modes}.

\begin{table}[h]
\centering
\caption{Approximate modes for selected values of \(n\), together with the
deficit \(\delta_n = e^{-e}-\mathrm{mode}\) from Euler's constant
\(e^{-e}\approx 0.0659880\).}
\label{tab:modes}
\footnotesize
\setlength{\tabcolsep}{3.5pt}
\begin{tabular}{lcccccccccc}
\toprule
\(n\) & 3 & 5 & 7 & 9 & 11 & 15 & 31 & 101 & 501 & 2001 \\
\midrule
Mode & 0.06399 & 0.06344 & 0.06316 & 0.06303 & 0.06298 & 0.06300 & 0.06346 & 0.06464 & 0.06560 & 0.06588 \\
\(10^{3}\delta_n\) & 2.00 & 2.55 & 2.83 & 2.96 & 3.01 & 2.98 & 2.53 & 1.35 & 0.38 & 0.11 \\
\bottomrule
\end{tabular}
\end{table}

Although these values are obtained numerically and therefore carry a small
discretization error, every one of them lies just below \(e^{-e}\approx 0.065988\).
The mode decreases until roughly height~11--15 and then rises again, approaching
\(e^{-e}\) from the left rather than drifting toward zero as the height
increases: the deficit \(e^{-e}-\text{mode}\) peaks near \(n=11\) and thereafter
falls by more than an order of magnitude between \(n=101\) and \(n=2001\).
This is not a numerical coincidence. Section~\ref{sec:limit} shows that
\(e^{-e}\) is exactly Euler's convergence boundary for infinite power towers and
that the limiting density has a one-sided singularity there, which forces the
modes to converge to \(e^{-e}\) from below.

\subsection{Mean and Variance}
\label{sec:meanvar}

Approximate means and variances, obtained by numerical integration, are
reported in Table~\ref{tab:meanvar}.

\begin{table}[h]
\centering
\caption{Approximate means and variances for selected values of \(n\). The final
column gives the corresponding quantities for the limit distribution
\(F_\infty\) of Section~\ref{sec:limit}.}
\label{tab:meanvar}
\footnotesize
\setlength{\tabcolsep}{3.5pt}
\begin{tabular}{lccccccccccc}
\toprule
\(n\) & 3 & 5 & 7 & 9 & 11 & 15 & 31 & 101 & 501 & 2001 & \(\infty\) \\
\midrule
Mean     & 0.4269 & 0.4024 & 0.3912 & 0.3850 & 0.3812 & 0.3769 & 0.3717 & 0.3692 & 0.3687 & 0.3687 & 0.3687 \\
Variance & 0.0911 & 0.0970 & 0.1007 & 0.1031 & 0.1048 & 0.1068 & 0.1096 & 0.1110 & 0.1113 & 0.1113 & 0.1113 \\
\bottomrule
\end{tabular}
\end{table}

As the height increases the mean declines slowly while the variance grows, both
quantities converging to the corresponding moments of the limit distribution
identified in Section~\ref{sec:limit}.

At first sight the behaviour of the variance is counter-intuitive: the mode
sharpens dramatically as \(n\) grows, and a sharper mode is normally associated
with \emph{less} dispersion, yet the variance \emph{increases} from \(0.0911\)
to \(0.1113\). The resolution is that the sharpening peak is not drawing mass in
from the whole interval; it is separating the distribution into two groups at
very different scales. As Theorem~\ref{thm:critical} below shows, in the limit a
proportion exactly \(1/e\approx 0.368\) of the total mass is compressed into the
short interval \([0,e^{-e}]\), a region of length under \(0.066\), while the
remaining \(1-1/e\approx 0.632\) is spread across \([e^{-e},1]\). It is this
growing separation between a tightly concentrated low-value group and a broadly
dispersed high-value group---rather than any spreading of the peak itself---that
drives the variance upward. The naive intuition ``sharper mode implies smaller
variance'' is reliable only for unimodal densities whose mass remains
concentrated about a single scale, which is precisely what fails here.

\subsection{Endpoint behaviour}

We now examine the density \(f_{n}(x)\) near the endpoints of the unit interval.

\paragraph{As \(x\to 0^+\).}
Although the power tower itself vanishes at the origin for every odd height
\(n\ge 3\), the density does not: evaluation of the recursion~\eqref{eq:densityrec}
shows that \(f_{n}(x)\to 1\) as \(x\to 0^{+}\) for every odd \(n\). From this
common boundary value the density rises extremely rapidly to its mode just below
\(e^{-e}\) and then falls. The taller the tower, the more violent this initial
spike: the peak becomes both taller and narrower as the height increases.

\paragraph{As \(x\to 1^-\).}
At the right endpoint the situation is far more regular. Direct evaluation of the recursive formula shows
\[
f_{n}(1) = 1
\]
for every height. Moreover the density approaches this value smoothly from below after the mid-interval dip. There is no singularity or rapid oscillation near \(1\); the behaviour is essentially the same for all odd heights once \(x\) is bounded away from zero.

\section{The Limit Distribution and Euler's Constant \texorpdfstring{\(e^{-e}\)}{e\textasciicircum(-e)}}
\label{sec:limit}

We now identify the limiting behaviour of the family as the height tends to
infinity through odd integers. This resolves the shape of the spike observed in
Section~\ref{sec:mode} and explains the role of the constant \(e^{-e}\).

\subsection{Euler's convergence boundary}

Fix \(x\in(0,1)\) and set \(t=\ln(1/x)>0\). The tower
recursion~\eqref{eq:towerdef} is the orbit of \(y_1=x\) under the map
\begin{equation}
    \varphi_x(y) = x^{\,y} = e^{-t y},
    \label{eq:phidef}
\end{equation}
so that \(y_m = \varphi_x^{\,m-1}(y_1)\). The map \(\varphi_x\) is a decreasing
contraction-or-expansion whose unique fixed point \(y^{*}(x)\) solves
\(y=x^{y}\), equivalently \(t y e^{t y} = t\), and is therefore given in closed
form by the Lambert \(W\) function~\cite{Corless1996} as
\begin{equation}
    y^{*}(x) = \frac{W(t)}{t}, \qquad t = \ln(1/x).
    \label{eq:fixedpoint}
\end{equation}
Since \(\varphi_x'(y) = -t\,e^{-ty}\), the multiplier at the fixed point is
\(\varphi_x'(y^{*}) = -t\,y^{*} = -W(t)\). The fixed point is therefore
attracting precisely when \(W(t)<1\), that is when \(t<e\), that is when
\(x>e^{-e}\); it loses stability exactly at
\begin{equation}
    x = e^{-e} \approx 0.0659880, \qquad t = e, \qquad y^{*} = \frac{W(e)}{e} = \frac{1}{e},
    \label{eq:critical}
\end{equation}
using \(W(e)=1\). This is Euler's classical convergence boundary for infinite
power towers~\cite{Euler1783,Knoebel1981}: the multiplier passes through \(-1\),
which is a period-doubling bifurcation. For \(x<e^{-e}\) the fixed point is
repelling and the orbit is attracted instead to a two-cycle
\(\{y_{-}(x),\,y_{+}(x)\}\) with \(y_{-}(x)<y^{*}(x)<y_{+}(x)\), the two points
being the solutions of \(y=\varphi_x(\varphi_x(y))\) other than \(y^{*}\).

Because \(y_1=x\) is small for small \(x\) while \(y_2=x^{x}\) is close to~\(1\),
the odd-index iterates track the \emph{lower} point of this two-cycle and the
even-index iterates track the upper one. This dichotomy is the exact analytic
counterpart of the even/odd distinction imposed in Section~\ref{sec:cons}.

The following observation is what makes the analysis elementary, and in
particular rules out any behaviour more complicated than a two-cycle.

\begin{lemma}
\label{lem:monotone}
For every \(x\in(0,1)\) the second iterate \(\psi_x := \varphi_x\circ\varphi_x\)
is strictly increasing on \((0,\infty)\). Consequently every orbit of
\(\psi_x\) is monotone and converges to a fixed point of \(\psi_x\).
Moreover \(y_3 > y_1\) for every \(x\in(0,1)\), so the odd-index subsequence
\(y_1,y_3,y_5,\dots\) is strictly increasing.
\end{lemma}

\begin{proof}
By the chain rule,
\[
\psi_x'(y) = \varphi_x'\bigl(\varphi_x(y)\bigr)\,\varphi_x'(y)
= \bigl(-t\,e^{-t\varphi_x(y)}\bigr)\bigl(-t\,e^{-ty}\bigr)
= t^{2}e^{-t\varphi_x(y)}e^{-ty} > 0 ,
\]
the product of two negative factors. A strictly increasing continuous self-map of
an interval has only monotone orbits: if \(\psi_x(y)>y\) then applying the
increasing map \(\psi_x\) repeatedly preserves the inequality, so the orbit is
increasing, and similarly in the decreasing case. Being monotone and bounded
(\(\varphi_x\) maps into \((0,1)\)), each orbit converges, necessarily to a fixed
point of \(\psi_x\) by continuity. This excludes cycles of period greater
than two for \(\varphi_x\), and in particular excludes any further
period-doubling cascade.

For the final claim, \(y_3>y_1\) reads \(x^{x^{x}}>x\). Taking logarithms gives
\(x^{x}\ln x > \ln x\), and dividing by \(\ln x<0\) reverses the inequality to
\(x^{x}<1\), which holds for every \(x\in(0,1)\).
\end{proof}

\begin{theorem}[Limit distribution]
\label{thm:limit}
Define \(F_\infty:[0,1]\to[0,1]\) by \(F_\infty(0)=0\), \(F_\infty(1)=1\) and,
for \(x\in(0,1)\) with \(t=\ln(1/x)\),
\[
F_{\infty}(x) =
\begin{cases}
y_{-}(x), & 0<x<e^{-e},\\[2pt]
\dfrac{W(t)}{t}, & e^{-e}\le x<1,
\end{cases}
\]
where \(y_{-}(x)\) denotes the lower point of the attracting two-cycle of
\(\varphi_x\). Then \(F_n \to F_\infty\) pointwise on \([0,1]\) as \(n\to\infty\)
through odd integers, \(F_\infty\) is a continuous cumulative distribution
function, the convergence is uniform on \([0,1]\), and consequently
\(X_n \Rightarrow X_\infty\).
\end{theorem}

\begin{proof}
\emph{Pointwise convergence.} Fix \(x\in(0,1)\). By Lemma~\ref{lem:monotone} the
odd-index subsequence \(y_1,y_3,y_5,\dots\) is a strictly increasing orbit of
\(\psi_x\) and therefore converges to the smallest fixed point of \(\psi_x\)
lying at or above \(y_1=x\).

We first record that
\begin{equation}
    x < y^{*}(x) \qquad\text{for every } x\in(0,1).
    \label{eq:xbelowfp}
\end{equation}
Indeed, writing \(w=W(t)>0\) so that \(t=we^{w}\), the inequality
\(x=e^{-t}<W(t)/t=y^{*}\) is equivalent to \(t\,e^{-t}<w\), that is to
\(e^{w}e^{-we^{w}}<1\), that is to \(w\bigl(1-e^{w}\bigr)<0\), which holds
because \(w>0\) and \(e^{w}>1\).

By Euler's classical convergence theorem~\cite{Euler1783,Knoebel1981}, the full
tower sequence \((y_m)_{m\ge1}\) converges exactly when \(x\ge e^{-e}\), and for
\(x<e^{-e}\) it splits into two subsequential limits forming a two-cycle of
\(\varphi_x\). We use this only to know \emph{whether} a two-cycle exists; the
identification of the branch is then forced by monotonicity, as follows.

If \(x\ge e^{-e}\), the whole sequence converges, necessarily to the unique
fixed point, so the odd subsequence has limit \(y^{*}(x)=W(t)/t\). (At
\(x=e^{-e}\) the fixed point is neutral and the convergence is algebraic rather
than geometric, but convergence still holds.)

If \(x<e^{-e}\), then \(\psi_x\) has the three fixed points
\(y_{-}<y^{*}<y_{+}\). Since \(\psi_x\) is increasing, the sign of
\(\psi_x(y)-y\) is constant between consecutive fixed points and alternates
across them; as \(\psi_x(0^{+})=e^{-t}>0\) we get \(\psi_x(y)>y\) on
\((0,y_{-})\) and \(\psi_x(y)<y\) on \((y_{-},y^{*})\). Orbits starting in
\((y_{-},y^{*})\) are therefore \emph{decreasing}. But by
Lemma~\ref{lem:monotone} our orbit is strictly increasing, and by
\eqref{eq:xbelowfp} it starts at \(x<y^{*}\); hence \(x\in(0,y_{-})\) and the
orbit increases to \(y_{-}(x)\).

\emph{\(F_\infty\) is a continuous CDF.} Each \(F_n\) is non-decreasing by
Theorem~\ref{thm:positivity}, so the pointwise limit \(F_\infty\) is
non-decreasing. As \(x\to 0^{+}\) we have \(t\to\infty\) and
\(0<F_\infty(x)\le y^{*}(x)=W(t)/t\to 0\), so \(F_\infty(0^{+})=0\); and
\(W(t)/t\to 1\) as \(t\to 0^{+}\), so \(F_\infty(1^{-})=1\). On \((e^{-e},1)\)
the function \(W(t)/t\) is real-analytic, and on \((0,e^{-e})\) the branch
\(y_{-}\) is real-analytic by the implicit function theorem, the relevant
derivative \(\partial_y(\psi_x(y)-y)\) being non-zero at \(y_{-}\) away from the
bifurcation point. Continuity at \(x=e^{-e}\) follows from the expansion of
Proposition~\ref{prop:sqrt} below, in which \(y_{-}(x)\to 1/e = y^{*}(e^{-e})\)
as \(x\uparrow e^{-e}\). Hence \(F_\infty\) is a continuous CDF on \([0,1]\).

\emph{Uniformity.} Pointwise convergence of distribution functions to a
continuous limit is exactly convergence in distribution, and by P\'olya's
theorem such convergence is automatically uniform on \([0,1]\).
\end{proof}

\begin{figure}[h!]
    \centering
    \includegraphics[width=0.78\textwidth]{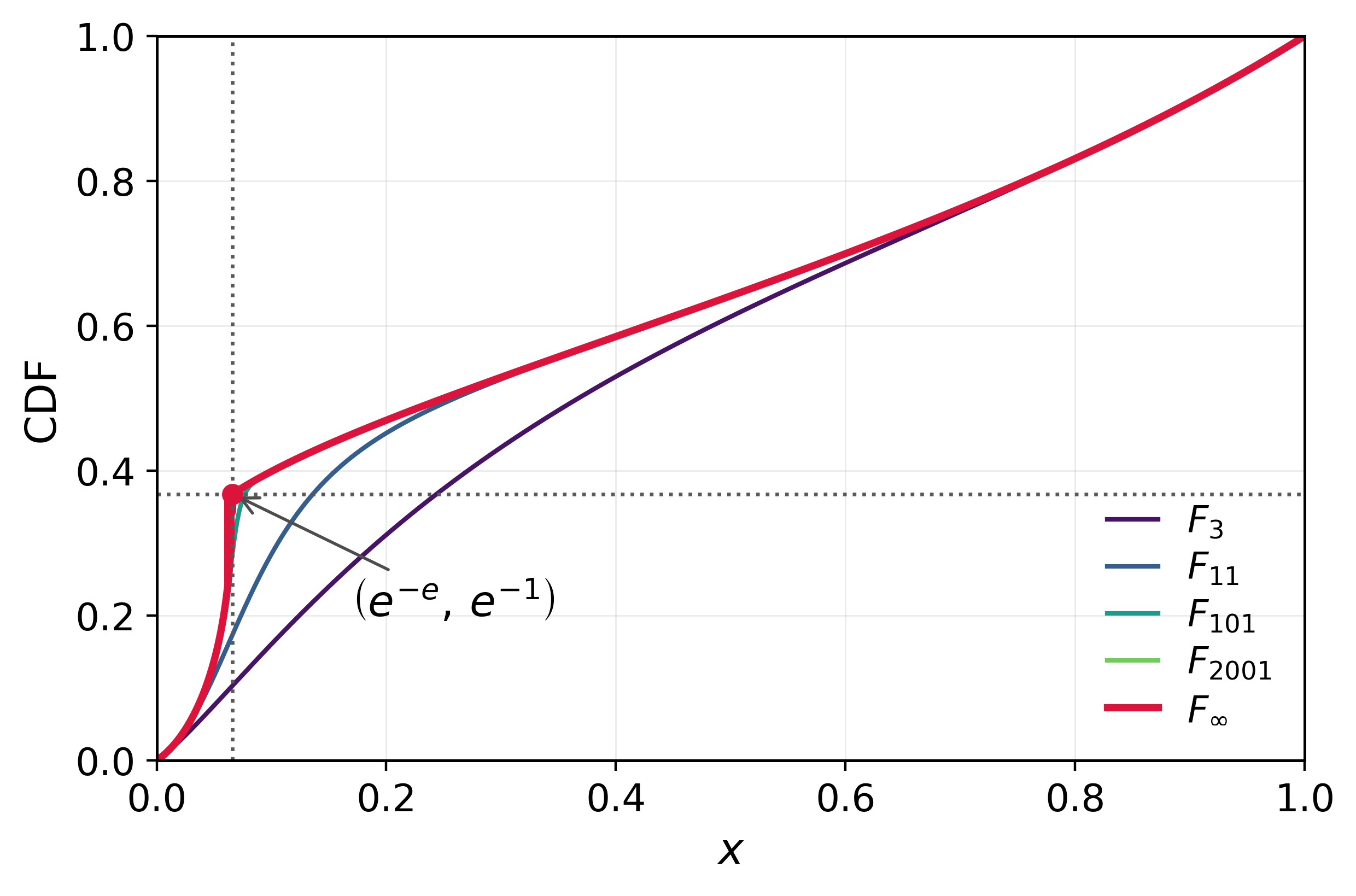}
    \caption{Odd-height tetration CDFs converging to the limit \(F_\infty\) of
    Theorem~\ref{thm:limit}. The dotted lines mark the critical point
    \(x=e^{-e}\) and the level \(e^{-1}\); by Theorem~\ref{thm:critical} the limit
    curve passes exactly through \((e^{-e},e^{-1})\), where it has a vertical
    tangent.}
    \label{fig:limit}
\end{figure}

Figure~\ref{fig:limit} displays the convergence. By height \(n=2001\) the tower
is numerically indistinguishable from \(F_\infty\) at the resolution of the plot.
Away from the critical point the convergence is geometric and very fast: in
double precision \(|F_{1001}-F_\infty|\) is already at the level of rounding
error (below \(10^{-15}\)) for \(x\le 0.06\) and for \(x\ge 0.1\). Inside that
window the rate degrades as \(e^{-e}\) is approached, the multiplier
\(|{-}W(t)|\) tending to~\(1\); at \(x=e^{-e}\) itself convergence is algebraic,
and \(|F_{1001}-F_\infty|\approx 2.8\times 10^{-2}\) there. This slow approach
precisely at the critical point is the finite-height signature of the
singularity described in Proposition~\ref{prop:sqrt}.

\subsection{The critical point}

\begin{theorem}[Mass below the critical point]
\label{thm:critical}
The limit distribution satisfies
\[
F_{\infty}\bigl(e^{-e}\bigr) = \frac{W(e)}{e} = \frac{1}{e},
\]
that is, \(X_\infty\) places mass exactly \(1/e\) on \([0,e^{-e}]\).
\end{theorem}

\begin{proof}
At \(x=e^{-e}\) we have \(t=e\), and \(W(e)=1\) since \(1\cdot e^{1}=e\).
Substituting into~\eqref{eq:fixedpoint} gives \(y^{*}=1/e\), and by
Theorem~\ref{thm:limit} this value is \(F_\infty(e^{-e})\).
\end{proof}

Theorem~\ref{thm:critical} is the exact statement underlying the appearance of
\(1/e\) in the numerical results of Section~\ref{sec:meanvar}. It is worth
emphasising that \(1/e\) arises here as a \emph{probability}, not as the mean;
see Remark~\ref{rem:notoneovere}.

\begin{proposition}[Square-root branching and the spike]
\label{prop:sqrt}
As \(x\uparrow e^{-e}\),
\[
F_{\infty}\bigl(e^{-e}-h\bigr) = \frac{1}{e} - C\sqrt{h} + O(h),
\qquad C = 2.1276555\ldots,
\]
so that the limiting density \(f_\infty = F_\infty'\) satisfies
\[
f_{\infty}(x) \sim \frac{C}{2}\,\bigl(e^{-e}-x\bigr)^{-1/2}
\qquad (x \uparrow e^{-e}),
\]
an integrable, one-sided inverse-square-root singularity. The density
\(f_\infty\) is finite and real-analytic elsewhere on \((0,1)\); in particular it
is bounded on \([e^{-e},1)\).
\end{proposition}

\begin{proof}[Proof sketch]
At the critical parameter the multiplier of \(\varphi_x\) at \(y^{*}\) equals
\(-1\), so the bifurcation is a generic period doubling: writing \(u=y-y^{*}\)
and \(t=e+\delta\), the second-iterate map \(\varphi_x\circ\varphi_x\) has the
normal form \(u \mapsto u + a\,\delta\,u + b\,u^{3} + O(\delta u^{2},u^{4})\) with
\(a,b\ne 0\), whence the two-cycle points separate from \(y^{*}\) at rate
\(|u|\asymp\sqrt{\delta}\). Since \(t=\ln(1/x)\) is analytic with
\(dt/dx = -1/x = -e^{e}\ne 0\) at \(x=e^{-e}\), we have
\(\delta\asymp h\) and therefore \(y^{*}-y_{-}\asymp\sqrt{h}\). Differentiating
the expansion gives the stated behaviour of \(f_\infty\). Analyticity away from
the critical point is the implicit function theorem, as in
Theorem~\ref{thm:limit}.
\end{proof}

The constant \(C\) was obtained numerically: computed in \(60\)-digit arithmetic,
the ratio \(\bigl(1/e-F_{\infty}(e^{-e}-h)\bigr)/\sqrt{h}\) increases
monotonically through the values \(2.127450\), \(2.127591\), \(2.127635\),
\(2.127649\) and \(2.127653\) at \(h=10^{-8},10^{-9},\dots,10^{-12}\). The
\(O(h)\) term in the expansion predicts that this ratio should approach \(C\)
from below at rate \(O(\sqrt{h})\), that is, by a factor of
\(\sqrt{10}\approx 3.16228\) per decade in \(h\). The observed successive
increments are \(1.4026\times10^{-4}\), \(4.4349\times10^{-5}\),
\(1.4024\times10^{-5}\) and \(4.4347\times10^{-6}\), shrinking by factors
\(3.16272\), \(3.16242\) and \(3.16232\)---converging to \(\sqrt{10}\), which
both confirms the square-root branching of Proposition~\ref{prop:sqrt} and
justifies Richardson extrapolation of the sequence. That extrapolation is itself
stable, returning \(2.127655467\), \(2.127655458\), \(2.127655457\) and
\(2.127655457\) from the four consecutive pairs, giving
\(C=2.1276555\ldots\). We have not identified a closed
form for \(C\); extracting one from the period-doubling normal form of
Proposition~\ref{prop:sqrt} is a natural next step.

Proposition~\ref{prop:sqrt} accounts for the entire qualitative picture of
Section~\ref{sec:mode}. The limiting density is unbounded as \(x\uparrow e^{-e}\)
and bounded on the other side, so the maximum of \(f_n\) must occur to the
\emph{left} of \(e^{-e}\) for every finite \(n\) and must migrate toward
\(e^{-e}\) as \(n\) grows. This is exactly the pattern of Table~\ref{tab:modes},
in which every mode lies below \(e^{-e}\) and the deficit decreases by more than
an order of magnitude from \(n=101\) to \(n=2001\). The finite-height densities
are smooth approximations to a function with an integrable singularity: the peak
grows without bound while the area it encloses tends to zero, which is why the
limit law has no atom at \(e^{-e}\) despite the diverging spike.

\subsection{Limiting moments}

\begin{proposition}[Limiting mean and variance]
\label{prop:moments}
With \(F_\infty\) as in Theorem~\ref{thm:limit},
\[
\mathbb{E}[X_{\infty}] = \int_{0}^{1}\bigl(1-F_{\infty}(x)\bigr)\,dx = 0.3686906\ldots,
\]
\[
\operatorname{Var}[X_{\infty}]
= \int_{0}^{1} 2x\bigl(1-F_{\infty}(x)\bigr)\,dx - \mathbb{E}[X_{\infty}]^{2}
= 0.1112922\ldots.
\]
\end{proposition}

\begin{proof}
The integral representations are the standard tail formulas for a non-negative
random variable bounded by~1, valid since \(F_\infty\) is a continuous CDF
supported on \([0,1]\) by Theorem~\ref{thm:limit}. The numerical values were
obtained by splitting the range at \(e^{-e}\) and applying the substitution
\(x = e^{-e}-u^{2}\) on the lower piece, which removes the square-root
singularity of Proposition~\ref{prop:sqrt} and renders both integrands smooth;
the quoted digits are stable under refinement and agree with direct quadrature
against the height-\(6401\) tower.
\end{proof}

\begin{remark}[The limiting mean is not \(1/e\)]
\label{rem:notoneovere}
The numerical proximity of \(\mathbb{E}[X_{\infty}]=0.3686906\ldots\) to
\(1/e=0.3678794\ldots\) is suggestive, and from the moderate heights of
Table~\ref{tab:meanvar} alone one might reasonably guess that the mean converges
to \(1/e\). It does not. The difference is \(8.1119\times 10^{-4}\), whereas our
three independent evaluations of \(\mathbb{E}[X_\infty]\)---adaptive quadrature
with the singularity split out, \(40\)-digit quadrature after the
desingularising substitution \(x=e^{-e}-u^{2}\), and direct integration against
the height-\(6401\) tower---agree with one another to within
\(8.2\times10^{-6}\); the discrepancy from \(1/e\) is therefore two orders of
magnitude larger than our numerical uncertainty. Moreover the finite
height means approach \(0.3686906\) rather than \(1/e\): the gap
\(\mathbb{E}[X_n]-1/e\) takes the values \(1.361\times10^{-3}\),
\(8.840\times10^{-4}\), \(8.205\times10^{-4}\) and \(8.124\times10^{-4}\) at
\(n=101,401,1601,6401\), plainly stabilising near \(8.1\times10^{-4}\) rather
than decaying to zero. The constant \(1/e\) does occur exactly in this problem,
but as the \emph{probability} \(F_{\infty}(e^{-e})\) of
Theorem~\ref{thm:critical}, and as the fixed point \(y^{*}=1/e\) at the
bifurcation~\eqref{eq:critical}, not as the limiting mean. We have not identified
closed forms for the two constants of Proposition~\ref{prop:moments}.
\end{remark}

Closed-form expressions for the moment-generating function and for moments of
order greater than two remain unknown, both for the finite heights and in the
limit.

\section{Statistical Inference}
\label{sec:inference}

The preceding sections treat the family analytically. We now turn to inference:
estimation of the height from data, a simulation study of the resulting
estimators, and an application to two real datasets bounded on the unit
interval. The outcome of the application is negative---the family is decisively
outperformed by the Beta and Kumaraswamy distributions---and
Section~\ref{sec:scope} explains structurally why, which we regard as the more
useful contribution of this section.

\subsection{Maximum likelihood}
\label{sec:mle}

Given a sample \(x_1,\dots,x_m\) from the unit interval, the log-likelihood of
the height-\(n\) tetration distribution is
\begin{equation}
    \ell(n) = \sum_{i=1}^{m} \log f_{n}(x_i),
    \label{eq:loglik}
\end{equation}
with \(f_n\) given by the recursion~\eqref{eq:densityrec}. The parameter is a
discrete one---an odd positive integer---so maximum likelihood is an exact
search over a grid of odd heights rather than a continuous optimisation, and no
regularity conditions or derivative computations are required.

A practical remark makes this cheap. Evaluating~\eqref{eq:densityrec} separately
for each candidate height costs \(O(N^{2})\) passes over the data when heights
up to \(N\) are considered, but the recursion builds every height on the way to
the tallest, so recording \((F_k,f_k)\) at each odd \(k\) during a single
sweep delivers the entire profile in \(O(N)\) passes. The whole profile over
\(n\le 101\) is therefore obtained at the cost of evaluating the tallest tower
once.

\paragraph{A two-parameter extension.}
The one-parameter family is very rigid (Section~\ref{sec:scope}), so we also
consider the exponentiated version, in the sense of Gupta and
Kundu~\cite{GuptaKundu2001},
\begin{equation}
    F_{n,\alpha}(x) = \bigl[\,x\uparrow\uparrow n\,\bigr]^{\alpha},
    \qquad \alpha>0,
    \label{eq:etdef}
\end{equation}
which is a valid CDF for every odd \(n\) and every \(\alpha>0\), and reduces to
the Beta\((\alpha,1)\) law at \(n=1\). Its density is
\(f_{n,\alpha}=\alpha F_n^{\alpha-1}f_n\), so
\[
\ell(n,\alpha) = m\log\alpha + (\alpha-1)\sum_{i}\log F_n(x_i)
                 + \sum_{i}\log f_n(x_i).
\]
For fixed \(n\) the score equation in \(\alpha\) solves in closed form,
\begin{equation}
    \hat\alpha(n) = \frac{-m}{\sum_{i=1}^{m}\log F_{n}(x_i)},
    \label{eq:alphahat}
\end{equation}
the denominator being strictly negative because \(0<F_n<1\) on \((0,1)\).
The two-parameter MLE is thus again an exact profile search over odd \(n\), with
\(\alpha\) eliminated analytically at each step; no numerical optimiser is
needed. Since \(\partial^{2}\ell/\partial\alpha^{2}=-m/\alpha^{2}\), the
observed information gives the conditional standard error
\(\widehat{\mathrm{se}}(\hat\alpha)=\hat\alpha/\sqrt{m}\), valid given
\(\hat n\).

\paragraph{Method of moments.}
Matching the sample mean to \(\mathbb{E}[X_n]\) is possible in principle but of
very limited use, because the attainable means form a short interval: by
Table~\ref{tab:meanvar} the mean decreases monotonically from \(0.5\) at
\(n=1\) to the limiting value \(0.3686906\ldots\) of
Proposition~\ref{prop:moments}, so
\begin{equation}
    \mathbb{E}[X_n] \in \bigl(0.36869,\; 0.5\bigr]
    \qquad\text{for every odd } n .
    \label{eq:meanrange}
\end{equation}
A moment estimator therefore exists only for samples whose mean lies in that
interval, which excludes both datasets analysed below. This is the first
symptom of the rigidity discussed in Section~\ref{sec:scope}. The exponentiated
family~\eqref{eq:etdef} does not suffer from it: its mean ranges over all of
\((0,1)\) as \(\alpha\) varies.

\subsection{Simulation}

Because \(F_n\) is continuous and strictly increasing on \((0,1)\) by
Theorem~\ref{thm:positivity}, inverse-transform sampling applies without
qualification: if \(U\sim\mathrm{Unif}(0,1)\) then \(F_n^{-1}(U)\) has the
height-\(n\) tetration distribution, and \(F_n^{-1}\) is evaluated by
bisection, which is unconditionally convergent here. Bisecting in \(\log x\)
rather than \(x\) is advisable, since for tall towers the lower quantiles are
extremely small. For the exponentiated family one simply draws
\(F_n^{-1}(U^{1/\alpha})\).

\subsection{Monte Carlo study}

Table~\ref{tab:mc} reports the behaviour of the maximum likelihood estimator of
the height, over \(2000\) replicates at each combination of true height and
sample size, with the search restricted to \(n\in\{1,3,\dots,41\}\). Bias and
mean squared error both decrease steadily with the sample size and the
probability of exact recovery increases, as consistency requires. Taller towers
are harder to identify, which is expected: the CDFs \(F_n\) converge to
\(F_\infty\), so neighbouring tall heights are increasingly similar. Even so the
profile likelihood is sharply peaked rather than flat---for a sample of \(5000\)
drawn at \(n=11\), the log-likelihood at \(n=15\) is already \(24\) units below
the maximum---so the height is genuinely identified, not merely bounded.

\begin{table}[h]
\centering
\caption{Monte Carlo performance of the maximum likelihood estimator
\(\hat n\) of the height, \(2000\) replicates per cell, search grid
\(n\in\{1,3,\dots,41\}\).}
\label{tab:mc}
\footnotesize
\setlength{\tabcolsep}{5pt}
\begin{tabular}{ccccccc}
\toprule
\(n\) & \(m\) & mean\((\hat n)\) & bias & MSE & \(\Pr(\hat n = n)\) \\
\midrule
3  & 20  & 4.993  & 1.993 & 47.562  & 0.391 \\
3  & 50  & 3.740  & 0.740 & 11.012  & 0.568 \\
3  & 100 & 3.332  & 0.332 & 1.452   & 0.718 \\
3  & 500 & 3.056  & 0.056 & 0.116   & 0.971 \\
\midrule
5  & 20  & 7.528  & 2.528 & 75.384  & 0.254 \\
5  & 50  & 6.009  & 1.009 & 18.410  & 0.351 \\
5  & 100 & 5.356  & 0.356 & 4.136   & 0.499 \\
5  & 500 & 5.160  & 0.160 & 0.620   & 0.848 \\
\midrule
11 & 20  & 13.730 & 2.730 & 116.668 & 0.115 \\
11 & 50  & 12.574 & 1.574 & 48.432  & 0.193 \\
11 & 100 & 11.738 & 0.738 & 20.464  & 0.222 \\
11 & 500 & 11.236 & 0.236 & 2.748   & 0.509 \\
\bottomrule
\end{tabular}
\end{table}

Table~\ref{tab:mcet} reports the corresponding study for the exponentiated
family, focusing on \(\hat\alpha\) from~\eqref{eq:alphahat}. The bias is small
already at \(m=50\) and the mean squared error falls at the expected rate.

\begin{table}[h]
\centering
\caption{Monte Carlo performance of \(\hat\alpha\) in the exponentiated
tetration family, \(2000\) replicates per cell.}
\label{tab:mcet}
\footnotesize
\setlength{\tabcolsep}{5pt}
\begin{tabular}{cccccccc}
\toprule
\(n\) & \(\alpha\) & \(m\) & mean\((\hat\alpha)\) & bias & MSE & \(\Pr(\hat n=n)\) \\
\midrule
3 & 1.5 & 20  & 1.6287 & 0.1287 & 0.2390 & 0.323 \\
3 & 1.5 & 50  & 1.5477 & 0.0477 & 0.1060 & 0.446 \\
3 & 1.5 & 100 & 1.5105 & 0.0105 & 0.0607 & 0.541 \\
3 & 1.5 & 500 & 1.5131 & 0.0131 & 0.0114 & 0.871 \\
\midrule
5 & 0.8 & 20  & 0.8543 & 0.0543 & 0.0575 & 0.208 \\
5 & 0.8 & 50  & 0.8180 & 0.0180 & 0.0206 & 0.298 \\
5 & 0.8 & 100 & 0.8126 & 0.0126 & 0.0097 & 0.429 \\
5 & 0.8 & 500 & 0.8037 & 0.0037 & 0.0018 & 0.806 \\
\bottomrule
\end{tabular}
\end{table}

\subsection{Application to data}
\label{sec:application}

We fit the family to two real datasets supported on the unit interval, both
distributed with the \texttt{datasets} package of \textsf{R}~\cite{RCore}.
The first is the pore-shape index of \(48\) rock samples from a petroleum
reservoir (\texttt{rock\$shape}; image analysis by R.\ Katz, University of
Oxford, from BP Research data), a dimensionless measure of pore cross-sectional
shape recorded in \((0,1)\).\footnote{The \textsf{R} documentation describes
this column as \(\text{perimeter}/\sqrt{\text{area}}\), but the stored values
are evidently normalised differently, lying in \((0.09,0.47)\) rather than above
\(2\sqrt{\pi}\). We therefore treat it simply as a dimensionless shape index on
the unit interval.} The second is the proportion of army draftees with education
beyond primary school in each of \(47\) French-speaking Swiss provinces around
1888 (\texttt{swiss\$Education}, divided by \(100\)), from Mosteller and
Tukey~\cite{MostellerTukey1977}. Both were chosen because they are
right-skewed with mass concentrated at small values, which is the regime in
which the tetration density is least unlike real data; no dataset on which the
family performs better was found.

\begin{figure}[h!]
    \centering
    \includegraphics[width=\textwidth]{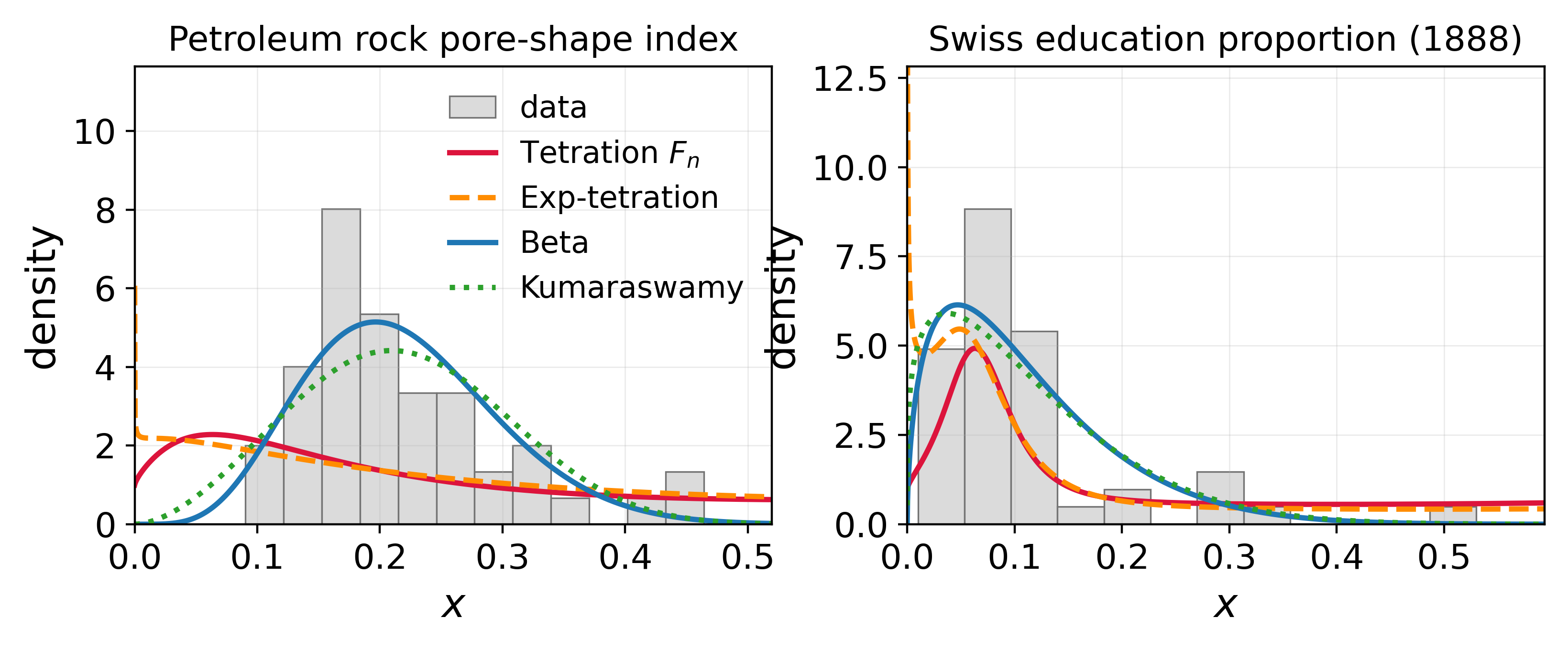}
    \caption{Fitted densities for the two datasets of Table~\ref{tab:app}.
    Left: the tetration mode is pinned near \(e^{-e}\approx 0.066\) while the
    data concentrate near \(0.18\). Right: the location is better matched, but
    the tetration right tail cannot decay, because \(f_n(1)=1\) for every
    height. The vertical range is truncated; several fitted densities diverge
    as \(x\to 0\).}
    \label{fig:appfit}
\end{figure}

We compare the one-parameter tetration family and its exponentiated extension
against the Beta distribution and the Kumaraswamy
distribution~\cite{Kumaraswamy1980}, the two standard two-parameter competitors
on the unit interval. Table~\ref{tab:app} reports maximised log-likelihoods,
AIC, BIC, and the Kolmogorov--Smirnov statistic with its \(p\)-value;
Figure~\ref{fig:appfit} shows the fitted densities.

\begin{table}[h]
\centering
\caption{Fits to two unit-interval datasets. Lower AIC and BIC are better;
\(k\) is the number of estimated parameters. The Kolmogorov--Smirnov
\(p\)-values reject both tetration models on both datasets and reject neither
competitor.}
\label{tab:app}
\footnotesize
\setlength{\tabcolsep}{4.5pt}
\begin{tabular}{llrrrrrl}
\toprule
Data & Model & \(\log L\) & \(k\) & AIC & BIC & KS & \(p\) \\
\midrule
\texttt{rock\$shape}      & Beta            & 55.603 & 2 & \(-107.205\) & \(-103.463\) & 0.1428 & 0.256 \\
\(m=48\)                  & Kumaraswamy     & 52.494 & 2 & \(-100.988\) & \(-97.246\)  & 0.1534 & 0.188 \\
                          & Tetration       & 12.891 & 1 & \(-23.783\)  & \(-21.912\)  & 0.4130 & \(<10^{-4}\) \\
                          & Exp-tetration   & 13.174 & 2 & \(-22.347\)  & \(-18.605\)  & 0.4099 & \(<10^{-4}\) \\
\midrule
\texttt{swiss\$Education} & Beta            & 59.607 & 2 & \(-115.214\) & \(-111.514\) & 0.1667 & 0.131 \\
\(m=47\)                  & Kumaraswamy     & 58.672 & 2 & \(-113.345\) & \(-109.644\) & 0.1734 & 0.105 \\
                          & Exp-tetration   & 45.780 & 2 & \(-87.559\)  & \(-83.859\)  & 0.3006 & 0.0003 \\
                          & Tetration       & 40.181 & 1 & \(-78.363\)  & \(-76.513\)  & 0.4385 & \(<10^{-4}\) \\
\bottomrule
\end{tabular}
\end{table}

The conclusion is unambiguous and negative. On both datasets the Beta and
Kumaraswamy fits are superior by \(40\)--\(80\) units of AIC, and the
Kolmogorov--Smirnov test rejects both tetration models at any conventional
level while rejecting neither competitor. The fitted heights are
\(\hat n = 5\) and \(\hat n = 19\); the exponentiated extension recovers only a
small part of the gap. We also fitted the remaining unit-interval datasets in
the same package, with the same outcome; on four of them the maximised
likelihood occurs at \(n=1\), which is the uniform distribution, so that the
family has nothing to offer beyond the trivial member. The one-parameter
tetration law is also beaten by one-parameter competitors such as the power law
\(F(x)=x^{a}\) and the truncated exponential, so the deficit is not an artefact
of comparing one parameter against two.

\subsection{Scope of applicability}
\label{sec:scope}

The failure is systematic rather than accidental, and its cause is visible in
Section~\ref{sec:dens}: the endpoint values of the density are the same for
every member of the family,
\[
f_n(0^{+}) = 1 \qquad\text{and}\qquad f_n(1) = 1 ,
\]
and the interior mode is confined to a neighbourhood of \(e^{-e}\), converging
to it as the height grows (Proposition~\ref{prop:sqrt}). The height \(n\)
controls only the sharpness of the spike. The family therefore has no freedom to
adjust either tail: it cannot produce a density vanishing at an endpoint, nor
one diverging there, whereas the Beta, Kumaraswamy and power families all do so
freely through their shape parameters. The right-hand panel of
Figure~\ref{fig:appfit} shows the consequence directly---the fitted tetration
density is still at height \(\approx 0.6\) at \(x=0.5\), where the data have
essentially no mass, purely because it must return to \(1\) at \(x=1\).

Two honest conclusions follow. First, the tetration family should not be
recommended as a general-purpose model for bounded data; on the evidence here
the Beta and Kumaraswamy families dominate it comprehensively. Second, the
family is plausibly useful only in the narrow regime its geometry actually
describes: data that are close to uniform over most of the interval but carry a
pronounced spike at small values near \(e^{-e}\), with no systematic decay
toward either endpoint. We have not identified a dataset of that kind. The
contribution of the present paper is accordingly the analytic structure of
Sections~\ref{sec:cons}--\ref{sec:limit}---in particular the limit law and the
role of Euler's constant---rather than a competitive modelling tool, and the
estimation theory above should be read as establishing what the family can and
cannot do.

\section{Qualitative Analogies}
\label{sec:apps}

There are a few natural phenomena that exhibit a qualitatively similar trend:
an initial response that overshoots a baseline level, drops rapidly, and then
gradually recovers. Two examples that first suggested themselves when examining
the form of the density are the classical illustration of the Dunning--Kruger
effect and the technology hype cycle.

In the usual schematic of the Dunning--Kruger effect the horizontal axis
represents an individual's actual knowledge of a subject while the vertical axis
represents expressed confidence. The resulting curve shows a high initial peak
at low knowledge, followed by a trough and a later rise. The tetration density
produces a qualitatively similar profile. Although more recent work has
questioned the empirical accuracy of that particular curve, the originally
postulated shape remains a useful morphological analogy for the density studied
here.

Likewise, the technology hype cycle traces a rapid rise in attention or
notoriety, a subsequent trough of disillusionment, and a gradual recovery.
Time series of public interest in novel technologies (for instance as measured
by Google Ngram frequencies) often display this same general shape. The
tetration density offers a simple one-parameter family that reproduces the
characteristic peak--dip--rise pattern.

Two caveats should be stated plainly. First, these analogies are morphological
only; no claim is made that either phenomenon is generated by the mechanism
studied here. Second, and more concretely, the curves just described are drawn on
arbitrary axes---time, competence, confidence---whereas the tetration
distribution lives strictly on \((0,1)\). Any use of the family on such a domain
therefore requires an affine transformation
\begin{equation}
    Y = a + (b-a)X, \qquad X\sim F_n, \qquad -\infty<a<b<\infty,
    \label{eq:affine}
\end{equation}
which has CDF \(F_Y(y)=F_n\bigl((y-a)/(b-a)\bigr)\) and density
\(f_Y(y) = (b-a)^{-1} f_n\bigl((y-a)/(b-a)\bigr)\) on \([a,b]\). Under this
transformation the mode sits at \(a+(b-a)\,\mathrm{mode}(f_n)\), so the
characteristic peak occurs at roughly \(6.6\%\) of the way along the range, and
the limiting mass below it is \(1/e\) by Theorem~\ref{thm:critical}. The
endpoints \(a\) and \(b\) are genuine additional parameters and would have to be
estimated alongside the height \(n\) in any applied use of the family. It should
be said plainly that the affine map rescales the interval but does not alter the
shape of the density, so it does not relieve the rigidity documented in
Section~\ref{sec:scope}; these analogies remain illustrations of a shape, not
evidence that the family models the phenomena.

\section{Conclusions}
\label{sec:conc}

We have introduced a new one-parameter family of distributions on the unit
interval generated by odd-height tetration. The cumulative distribution function
is given explicitly by the power tower of height \(n\), and the density satisfies
a simple recursive relation. Numerical investigation reveals a characteristic
shape consisting of a sharp spike near the origin, a subsequent dip, and
a gradual rise toward the right endpoint.

The main structural result is the identification of the limit as the height tends
to infinity through odd integers. The limiting CDF is \(W(t)/t\) with
\(t=\ln(1/x)\) above Euler's constant \(e^{-e}\) and the lower branch of the
attracting two-cycle of \(y\mapsto x^{y}\) below it, and the convergence is
uniform. The constant \(e^{-e}\) is precisely the parameter value at which
\(y\mapsto x^{y}\) undergoes a period-doubling bifurcation, which produces a
one-sided inverse-square-root singularity in the limiting density; this explains
the spike observed at every finite height and shows that the modes converge to
\(e^{-e}\) from the left. At the critical point the limit law satisfies
\(F_\infty(e^{-e})=1/e\) exactly.

The limiting mean and variance are \(0.3686906\ldots\) and \(0.1112922\ldots\).
The mean is close to \(1/e\) but is provably different from it; the constant
\(1/e\) enters this problem as a probability rather than as a mean. No closed
form has been identified for either limiting moment, nor for the branching
constant \(C\) of Proposition~\ref{prop:sqrt}, and closed-form expressions for
the moment-generating function and for moments of order greater than two remain
unknown.

On the inferential side, Section~\ref{sec:inference} shows that estimation
itself poses no difficulty. The height is discrete, so maximum likelihood is an
exact profile search, computable for all heights up to \(N\) in \(O(N)\) passes
over the data; in the exponentiated extension the second parameter is removed in
closed form; and a Monte Carlo study confirms that bias and mean squared error
decay as consistency requires. What fails is the fit. On two real unit-interval
datasets the Beta and Kumaraswamy distributions outperform the family by
\(40\)--\(80\) units of AIC, and a Kolmogorov--Smirnov test rejects the
tetration models while rejecting neither competitor. The reason is structural:
\(f_n(0^{+})=f_n(1)=1\) for every height, so no member of the family can adapt
either tail, and the height governs only the sharpness of the spike near
\(e^{-e}\). We therefore do not recommend the family as a general-purpose model
for bounded data, and we have stated in Section~\ref{sec:scope} the narrow
geometry for which it could be suitable.

Several analytic directions remain open. Closed forms for the two limiting
moments and for the branching constant \(C\) of Proposition~\ref{prop:sqrt}, and
the rate at which the finite-height modes approach \(e^{-e}\), are natural
questions suggested by the period-doubling picture. A further possibility,
suggested by the rigidity just described, is to enlarge the family by composing
it with a flexible transformation of the unit interval---the affine
map~\eqref{eq:affine} being the simplest---so as to free the endpoint behaviour
while retaining the tetration structure; whether anything of the characteristic
geometry survives such an enlargement is not clear to us.

The distinctive peak--dip--rise profile of the density bears a qualitative
resemblance to certain temporal or competence-related curves that exhibit an
early overshoot followed by recovery. These morphological parallels are noted
only as illustrations of the shape; formal statistical modelling of such
phenomena is left for future work.

\section*{Declarations}

\bmhead{Funding}
The author declares that no funds, grants, or other support were received during the preparation of this manuscript.

\bmhead{Conflict of interest}
The author has no relevant financial or non-financial interests to disclose.

\bmhead{Data availability}
The two datasets analysed in Section~\ref{sec:application} are publicly available and require no special access: both are distributed with the \texttt{datasets} package of \textsf{R}~\cite{RCore} as \texttt{rock\$shape} and \texttt{swiss\$Education}, the latter originating with Mosteller and Tukey~\cite{MostellerTukey1977}. All other numerical results are reproducible directly from the formulas given herein; the scripts used to generate the figures, tables and simulation study are available from the author upon reasonable request.

\bibliography{refs}

\end{document}